\documentclass{amsart}

\usepackage{amsmath,amssymb,amscd,amsfonts}

\newtheorem{theorem}{Theorem}
\newcommand{\bt}{\begin{theorem}}
\newcommand{\et}{\end{theorem}}
\newtheorem{lemma}{Lemma}
\newcommand{\bl}{\begin{lemma}}
\newcommand{\el}{\end{lemma}}
\newtheorem{corollary}{Corollary}
\newcommand{\bc}{\begin{corollary}}
\newcommand{\ec}{\end{corollary}}

\newcommand{\beq}{\begin{equation}}
\newcommand{\eeq}{\end{equation}}
\newcommand{\benum}{\begin{enumerate}}
\newcommand{\eenum}{\end{enumerate}}
\newcommand{\N}{\ensuremath{ \mathbf N }}

\newcommand{\mba}{\ensuremath{ \mathbf a}}
\newcommand{\mbb}{\ensuremath{ \mathbf b}}

\newcommand{\mbx}{\ensuremath{ \mathbf x}}
\newcommand{\mby}{\ensuremath{ \mathbf y}}

\newcommand{\bsmallmat}{\left(\begin{smallmatrix}}
\newcommand{\esmallmat}{\end{smallmatrix}\right)}

\newcommand{\bmat}{\left(\begin{matrix}}
\newcommand{\emat}{\end{matrix}\right)}

\DeclareMathOperator{\qqand}{\qquad\text{and}\qquad}

\title{Landau's converse to  H\" older's inequality}
\author{Melvyn B.   Nathanson} 
\address{Lehman College (CUNY), Bronx, NY 10468}
\email{melvyn.nathanson@lehman.cuny.edu}

\date{\today}



\begin{document}

\begin{abstract} 
By H\" older's inequality, 
if $\mathbf{x} \in  \ell^p$, then  $\mathbf{x}\mathbf{y} \in \ell^1$ for all $\mathbf{y} \in \ell^q$.  
Landau proved the converse result:  If $\mathbf{x}\mathbf{y} \in \ell^1$ for all $\mathbf{y} \in \ell^q$, 
then  $\mathbf{x} \in \ell^p$.  
This paper proves Landau's theorem and considers a related problem for Hilbert's inequality.  
\end{abstract}

\maketitle

\section{H\" older's inequality}
In this paper we prove H\" older's inequality and Landau's converse 
to H\" older's inequality.\footnote{These are the notes of a talk in the New York Number 
Theory Seminar on March 28, 2024.}  
There is a related inequality of Hilbert (proved in Appendix~\ref{Landau:appendix:Hilbert} ).  
The existence of a converse to Hilbert's inequality that is similar to Landau's theorem  
seems to be an open problem.  

The classic reference to inequalities is Hardy, Littlewood, and P\" olya~\cite{hard-litt-poly88}
A recent introduction is Steele~\cite{stee13}. 

There are many proofs of H\" older's inequality.  
Here is a simple proof of the finite form of the inequality. 

\bt                                \label{Landau:theorem:finite}
Let $p>1$ and $q>1$ be real numbers such that 
\[
\frac{1}{p}+\frac{1}{q}=1.
\]
For all positive real numbers $x_1,\ldots, x_n,y_1,\ldots, y_n$, 
\[
\sum_{i=1}^n x_iy_i \leq 
\left( \sum_{i=1}^n x_i^p \right)^{1/p} \left( \sum_{i=1}^n y_i^q \right)^{1/q}
\]
\et

\begin{proof}
The function $f(x)$ is \emph{convex} on the interval $(u,v)$ if 
\[
f(ta+(1-t)b) \leq tf(a) + (1-t)f(b)
\]
for all $a,b$ such that $u < a < b < v$.
The function  $f(x)$ is \emph{concave} on the interval $(u,v)$ 
if $-f(x)$ is convex on $(u,v)$ or, equivalently, if 
\[
f(ta+(1-t)b) \geq tf(a) + (1-t)f(b)
\]
for all $a,b$ such that $u < a < b < v$.
The function $f(x) = \log x$ is concave on $(0,\infty)$.  
For all $i \in \{1,2,\ldots, n\}$, choosing  $a = x_i^p$,   $b = y_i^q$,  
$t = 1/p$, and $1-t = 1/q$, 
we obtain 
\[
\log\left( \frac{1}{p} x_i^p + \frac{1}{q} y_i^q \right)  
\geq  \frac{1}{p} \log x_i^p +  \frac{1}{q} \log y_i^q = \log x_iy_i. 
\]
Exponentiating this inequality gives a simple form of \emph{Young's inequality}: 
\[
x_iy_i \leq  \frac{1}{p} x_i^p + \frac{1}{q} y_i^q. 
\]
Summing over $i$ we obtain 
\beq                 \label{Landau:ineq-xypq}
\sum_{i=1}^n x_iy_i \leq  \frac{1}{p} \sum_{i=1}^n x_i^p + \sum_{i=1}^n \frac{1}{q} y_i^q. 
\eeq
Let 
\[
X = \left( \sum_{i=1}^n x_i^p \right)^{1/p} \qqand Y = \left( \sum_{i=1}^n y_i^q \right)^{1/q}.
\]
Replacing $x_i$ with $x_i/X$ and $y_i$ with $y_i/Y$ in~\eqref{Landau:ineq-xypq}, we obtain 
\[
\frac{1}{XY} \sum_{i=1}^n x_iy_i \leq  
\frac{1}{p} \frac{\sum_{i=1}^n x_i^p}{X^p} + \frac{1}{q} \sum_{i=1}^n \frac{\sum_{i=1}^n y_i^q}{Y^q}  
= \frac{1}{p}  + \frac{1}{q} = 1.
\]
This completes the proof.  
\end{proof}

For sequences $\mbx = (x_i)_{i=1}^{\infty}$ and $\mby = (y_i)_{i=1}^{\infty}$ 
of real or complex numbers, 
we define the product sequence $\mbx \mby  = (x_i y_i)_{i=1}^{\infty}$.  
For all $p \geq 1$, the $p$-norm of the sequence $\mbx = (x_i)_{i=1}^{\infty}$  is 
\[
\| \mbx \|_p = \left( \sum_{i=1}^{\infty} |x_i|^p \right)^{1/p}.
\]
The Lebesgue space $\ell^p$ is the set of all sequences \mbx\ such that $\|\mbx\|_p < \infty$. 

Theorem~\ref{Landau:theorem:finite} immediately implies the infinite form of H\" older's inequality. 

Let $\N = \{1,2,3,\ldots\}$ be the set of positive integers. 

\bt                                                                  \label{Landau:theorem:infinite}
Let $p>1$ and $q>1$ be real numbers such that 
\[
\frac{1}{p}+\frac{1}{q}=1.
\]
Let $\mbx = (x_i)_{i=1}^{\infty}$ and $\mby = (y_i)_{i=1}^{\infty}$ be sequences 
of real or complex numbers.   
If $\mbx = (x_i)_{i=1}^{\infty} \in \ell^p$ and $\mby = (y_i)_{i=1}^{\infty} \in \ell^q$, 
then 
\[
\|\mbx\mby\|_1 \leq \|\mbx \|_p \|\mby\|_q.
\]
\et

\begin{proof}
Without loss of generality, we can assume that $x_i \neq 0$ and $y_i \neq 0$ for all $i \in \N$. 
By Theorem~\ref{Landau:theorem:finite}, for all $n \in \N$ we have 
\begin{align*} 
\sum_{i=1}^n |x_iy_i|&  \leq \left( \sum_{i=1}^n |x_i|^p \right)^{1/p} \left( \sum_{i=1}^n |y_i|^q \right)^{1/q} \\
& \leq \left( \sum_{i=1}^{\infty} |x_i|^p \right)^{1/p} \left( \sum_{i=1}^{\infty} |y_i|^q \right)^{1/q} \\
& =  \|\mbx \|_p \|\mby\|_q 
\end{align*} 
and so 
\[
\|\mbx\mby\|_1 = \sum_{i=1}^{\infty} |x_iy_i| \leq   \|\mbx \|_p \|\mby\|_q.
\]
This completes the proof.
\end{proof}

For all $p \geq 1$, the Lebesgue space  $\ell^p$ contains every finite sequence  and so  Theorem~\ref{Landau:theorem:infinite}, applied to the sequences 
$\mbx = (x_i)_{i=1}^n$ and $\mby = (y_i)_{i=1}^n$,
 gives Theorem~\ref{Landau:theorem:finite}.  
 Thus, Theorems~\ref{Landau:theorem:finite} and~\ref{Landau:theorem:infinite} are equivalent.

\section{Landau's theorem}

The Cauchy-Schwarz inequality states that if $\mbx \in \ell^2$ and $\mby \in  \ell^2$, then 
$\mbx \mby \in \ell^1$ and $\| \mbx \mby\|_1 \leq \| \mbx\|_2 \| \mby\|_2$.  
This is the case $p=q=2$ of H\" older's inequality.  
In 1906, motivated by Hilbert's papers on integral equations, 
Hellinger and T\" oplitz~\cite{hell-topl06} stated (but did not prove)  
the following converse of the Cauchy-Schwarz inequality.

\bt[Hellinger-T\" oplitz]
If $\mbx$ is a sequence   such that 
$\mbx \mby \in \ell^1$  for all sequences    $\mby \in  \ell^2$, 
then  $\mbx \in  \ell^2$.
\et
  
In 1907, Landau~\cite{land07} generalized this result  
and proved a  converse to the infinite H\" older's inequality.    
He wrote:
\begin{quotation}
Indem ich mir einen Beweis f\" ur diesen Satz zurechtlegte, fand ich folgende Verallgemeinerung, welche vielleicht f\" ur die Zwecke der Herren Hellinger und T\" oplitz n\" utzlich sein kann, und in welcher der oblige Hilfssatz als Specialfall $p=2$ enthalten ist.
\end{quotation}

Riesz included Landau's theorem in his book~\cite[pp.46--48]{ries13} 
on infinite systems of linear equations in infinitely many variables.  
For a history of H\" older's inequality, see Maligranda~\cite {mali98} 
and for recent related results, see Albuquerque, et al.\cite{albu17}.

Landau's proof uses the  following convergence-divergence theorem.

\bt [Abel-Dini]                                                           \label{Landau:theorem:Dini-1}
Let $\sum_{n=1}^{\infty} x_n$ be a divergent series of positive real numbers  
 and let $s_n = \sum_{i=1}^n x_i$ be the $n$th partial sum 
of the series. The series 
\beq                                                                      \label{Landau:Dini-1}
\sum_{n=1}^{\infty} \frac{x_n}{s_n^{\alpha}}  
\eeq
 diverges  if $\alpha \leq 1$ and converges if $\alpha > 1$.
\et

\begin{proof}
 Appendix~\ref{Landau:appendix:Dini}.
\end{proof}

\bt[Landau]

Let $p>1$ and $q>1$ be real numbers such that 
\[ 
\frac{1}{p} + \frac{1}{q} = 1.
\]
If $\mbx$ is a sequence  such that 
$\mbx \mby \in \ell^1$  for all sequences    $\mby \in  \ell^q$, 
then  $\mbx \in  \ell^p$.
\et

\begin{proof}
Without loss of generality, we can assume that 
 $\mbx = (x_i)_{i=1}^{\infty}$  
is a sequence of positive real numbers.  
We shall prove that if $\mbx \notin \ell^p$, then there exists $\mby \in \ell^q$ such that 
$\mbx\mby \notin \ell^1$.  

If $\mbx \notin \ell^p$, then the infinite series 
\[
\sum_{i=1}^{\infty} x_i^p 
\] 
diverges.  
For all $n \in \N$, the $n$th partial sum of the series is 
\[
s_n = \sum_{i=1}^n x_i^p.  
\] 
 Let  $s_0 = 0$.  
Then 
\[
x_i^p = s_i - s_{i-1} 
\]
for all $i \in \N$ and 
\[
0 = s_0 <  s_1 <  \cdots <  s_n <  s_{n+1} <  \cdots. 
\]
The divergence of the series $\sum_{i=1}^{\infty} x_i^p$ implies  
\[
\lim_{n\rightarrow \infty} s_n = \infty.
\]
By the Abel-Dini theorem, the series 
\beq                \label{Landau:Abel-Dini-div}
\sum_{i=1}^{\infty} \frac{x_i^p}{s_i} 
\eeq
diverges and the series 
\beq                \label{Landau:Abel-Dini-conv}
\sum_{i=1}^{\infty} \frac{x_i^p}{s_i^q} 
\eeq 
converges. Define $\lambda > 0$ by 
\[
 \frac{1}{\lambda^q} = \sum_{i=1}^{\infty} \frac{x_i^p}{s_i^q}.
\]

For all $i \in \N$, let $y_i$ be the positive real number 
\[
y_i = \lambda \frac{ \sqrt[q]{x_i^p}}{s_i}. 
\] 
and let $\mby   = (y_i)_{i=1}^{\infty}$.  
We have 
\[
\sum_{i=1}^{\infty} y_i^q = \lambda^q \sum_{i=1}^{\infty} \frac{x_i^p}{s_i^q} = 1
\]
and so $\mby = (y_i)_{i=1}^{\infty} \in \ell^q$.  
Because 
\[
1 + \frac{p}{q} = p 
\]
we have 
\begin{align*}
\sum_{i=1}^{\infty} x_iy_i  
&= \lambda \sum_{i=1}^{\infty} x_i \frac{ \sqrt[q]{x_i^p}}{s_i} \\ 
& = \lambda \sum_{i=1}^{\infty} \frac{x_i^{(1+ p/q)}}{s_i} \\
& = \lambda \sum_{i=1}^{\infty} \frac{x_i^{p}}{s_i} 
\end{align*}
This series diverges by Abel-Dini~\eqref{Landau:Abel-Dini-div}. 
Thus, if $\mbx \notin \ell^p$, then there exists $\mby \in \ell^q$ 
such that $\mbx\mby \notin \ell^1$. 
Equivalently, if $\mbx\mby \in \ell^1$ for all $\mby \in \ell^q$, then $\mbx \in \ell^p$. 
This completes the proof. 
\end{proof}

\section{Converse to the finite H\" older inequality }
Here is a converse to Theorem~\ref{Landau:theorem:finite}.  

\bt                                \label{Landau:theorem:finite-converse}
Let $p>1$ and $q>1$ be real numbers such that 
\[
\frac{1}{p}+\frac{1}{q}=1.
\]
If $x_1,\ldots, x_n,$ and $C$ are positive real numbers 
such that  
\[
\sum_{i=1}^n x_iy_i \leq C \left( \sum_{i=1}^n y_i^q \right)^{1/q}
\] 
for all positive real numbers $y_1,\ldots, y_n$, then 
\[
C \geq \left( \sum_{i=1}^n x_i^p \right)^{1/p}.
\]
\et

\begin{proof}

For all $i \in \{1,2,\ldots, n\}$, let $y_i$ be the   positive real number 
such that $x_iy_i = y_i^q$ or, equivalently,  
\[
y_i = x_i^{1/(q-1)}.   
\]
Because 
\[
p = \frac{q}{q-1} 
\]
we have  
\[
\sum_{i=1}^n x_iy_i = \sum_{i=1}^n x_i x_i^{1/(q-1)}  = \sum_{i=1}^n  x_i^{q/(q-1)} 
= \sum_{i=1}^n  x_i^p 
\]
and 
\[
 \left( \sum_{i=1}^n y_i^q \right)^{1/q} 
 =  \left( \sum_{i=1}^n   x_i^{q/(q-1)}   \right)^{1/q}
 =  \left( \sum_{i=1}^n   x_i^p  \right)^{1/q}. 
 \]
 Therefore, 
 \[
 \sum_{i=1}^n  x_i^p = \sum_{i=1}^n x_iy_i \leq C  \left( \sum_{i=1}^n y_i^q \right)^{1/q} 
 = C \left( \sum_{i=1}^n   x_i^p  \right)^{1/q}
 \]
 and 
 \[
  \left( \sum_{i=1}^n   x_i^p  \right)^{1/q} =   \left( \sum_{i=1}^n   x_i^p  \right)^{1-1/p}  \leq C. 
 \]
This completes the proof. 
\end{proof}

\appendix

\section{Dini's convergence theorem}       \label{Landau:appendix:Dini} 
Landau's converse to the infinite  H\" older inequality uses a convergence-divergence 
result   called the Abel-Dini theorem.  
It was proved by Dini in 1867 and is called the Abel-Dini theorem 
because of a related result proved by Abel in 1828.  

\bt                                                        Let $\sum_{n=1}^{\infty} a_n$ be a divergent series of positive terms 
 and let $s_n = \sum_{i=1}^n a_i$ be the $n$th partial sum 
of the series. The series 
\beq                                                                      \label{Landau:Dini-2}
\sum_{n=1}^{\infty} \frac{a_n}{s_n^{\alpha}}  
\eeq
 diverges  if $\alpha \leq 1$ and converges if $\alpha > 1$.
\et

\begin{proof}
The proof is from Knopp~\cite{knop56}.
Because the infinite series  $\sum_{n=1}^{\infty} a_n$ diverges, 
the sequence of partial sums $(s_n)_{n=1}^{\infty}$ is strictly increasing and  
$\lim_{n\rightarrow \infty} s_n = \infty$.  
It follows that for every positive integer $k$ there exists a positive integer $\ell$ such that 
$s_{k+\ell} > 2s_k$ or, equivalently, 
\[
\frac{s_k}{s_{k+\ell}} < \frac{1}{2}. 
\]
For all $s_k \geq 1$ and $\alpha \leq 1$ we have 
\[
\sum_{i=k+1}^{k+\ell}\frac{a_i}{s_i^{\alpha} }  \geq \sum_{i=k+1}^{k+\ell}\frac{a_i}{s_i} 
 \geq \sum_{i=k+1}^{k+\ell}\frac{a_i}{s_{k+\ell}} 
=  \frac{s_{k+\ell} - s_k}{s_{k+\ell}} = 1 - \frac{ s_k}{s_{k+\ell}} > \frac{1}{2}.  
\] 
and so the series~\eqref{Landau:Dini-2} diverges.  

Next we consider the case  $\delta > 0$ and $\alpha = 1 + \delta$. 
For all $r \in \N$ and $x \in (0,1)$, we have 
\[
1-x^r = (1-x) \left(1+x+x^2 + \cdots + x^{r-1} \right) < r(1-x).  
\] 
Choose $r \in \N$ such that $0 < 1/r < \delta$.  For all $n \geq 2$ we have 
\[
x = \left( \frac{s_{n-1}}{s_n} \right)^{1/r}  \in (0,1) 
\]
and    
\[
\frac{a_n}{s_n} = \frac{s_n - s_{n-1}}{s_n} 
= 1- \frac{s_{n-1}}{s_n} < r  \left( 1-  \left( \frac{s_{n-1}}{s_n} \right)^{1/r}\right).  
\]
We obtain the telescoping sum 
\begin{align*}
\sum_{n=1}^{\infty} \frac{a_n}{s_n^{\alpha}}  
& = \sum_{n=1}^{\infty} \frac{a_n}{s_n s_n^{\delta}}   \\ 
& \leq  \frac{1}{a_1^{1/r} } + \sum_{n=2}^{\infty} \frac{1}{ s_{n-1}^{1/r}} \frac{a_n}{s_n} \\
& =  \frac{1}{a_1^{1/r} } 
+ r \sum_{n=2}^{\infty}  \frac{1}{ s_{n-1}^{1/r}}   \left( 1-  \left( \frac{s_{n-1}}{s_n} \right)^{1/r}\right) \\ 
& =  \frac{1}{a_1^{1/r} } 
+ r \sum_{n=2}^{\infty}  \left(  \frac{1}{ s_{n-1}^{1/r}}  -  \frac{1}{ s_n^{1/r}} \right) \\
& =  \frac{r+1}{a_1^{1/r} }
\end{align*} 
and so the series~\eqref{Landau:Dini-1} converges.  
This completes the proof.  
\end{proof}

\section{Hilbert's inequality}         \label{Landau:appendix:Hilbert}

\bt[Hilbert] 
For all $\mba = (a_m)_{m=1}^{\infty} \in \ell^2$ and $\mbb = (b_n)_{n=1}^{\infty} \in \ell^2$,  
\[
\sum_{m=1}^{\infty} \sum_{n=1}^{\infty} \frac{a_m b_n}{m+n} < \pi \left(\sum_{m=1}^{\infty} a_m^2 \right)^{1/2}  \left(\sum_{n=1}^{\infty} b_n^2 \right)^{1/2}. 
\]
The constant $\pi$ is best possible. 
\et

In Hilbert's original proof, the constant on the right side of the inequality was $2\pi$.  
Schur obtained the constant $\pi$ and proved that this is  best possible.

There is a vast literature on Hilbert's inequality and its generalizations.  
The following $pq$ result is sometimes called the Hardy-Hilbert inequality.  
Here is a short proof that uses only H\" older's inequality.   

\bt
Let $p >1$ and $q>1$ satisfy 
\[
\frac{1}{p} + \frac{1}{q} = 1. 
\]
For all $\mba = (a_m)_{m=1}^{\infty} \in \ell^p$ and $\mbb = (b_n)_{n=1}^{\infty} \in \ell^q$,  
\[
\sum_{m=1}^{\infty} \sum_{n=1}^{\infty} \frac{a_m b_n}{m+n} < \frac{\pi}{\sin(\pi/p)} \left(\sum_{m=1}^{\infty} a_m^p \right)^{1/p}  \left(\sum_{n=1}^{\infty} b_n^q \right)^{1/q}.
\]
The constant $\pi/\sin(\pi/p)$ is best possible. 
\et

\begin{proof} 
We begin with two identities.  First, 
\[
\sin  \left(\frac{\pi}{p}\right) = \sin \pi \left(1-\frac{1}{p}\right) = \sin  \left(\frac{\pi}{q}\right) .
\]
Second, for $0 < \lambda < 1/p$, 
\[
\int_0^{\infty}  \frac{1}{ (1+t) t^{p \lambda}}  dt  = \frac{\pi}{\sin\left(\pi p\lambda \right)}. 
\]
Borwein~\cite{borw08} gives two proofs of this formula.  
\emph{Maple} and \emph{Mathemtica} also give this result.

For $\lambda > 0$, let 
\[
x_m = \frac{a_m}{(m+n)^{1/p} } \left( \frac{m}{n}\right)^{\lambda}
\qqand 
y_n = \frac{b_n}{ (n+m)^{1/q}  } \left( \frac{n}{m}\right)^{\lambda}.
\]
Then 
\begin{align*}
 \sum_{m=1}^{\infty} \sum_{n=1}^{\infty} \frac{a_m b_n}{m+n}  
 & = 
 \sum_{m=1}^{\infty} \sum_{n=1}^{\infty} \frac{a_m}{(m+n)^{1/p} } \left( \frac{m}{n}\right)^{\lambda}
\frac{b_n}{(m+n)^{1/q} }  \left( \frac{n}{m}\right)^{\lambda} \\
& = \sum_{m=1}^{\infty} \sum_{n=1}^{\infty} x_m y_n \\ 
& \leq \left( \sum_{m=1}^{\infty}   \sum_{n=1}^{\infty} x_m^p \right)^{1/p} 
 \left( \sum_{n=1}^{\infty} \sum_{m=1}^{\infty} y_n^q \right)^{1/q} 
 \hspace{0.3cm} \text{(by H\" older's inequality)}       \\
& = \left( \sum_{m=1}^{\infty}   \sum_{n=1}^{\infty} 
 \frac{a_m^p}{m+n } \left( \frac{m}{n}\right)^{p \lambda} \right)^{1/p}
 \left( \sum_{n=1}^{\infty} \sum_{m =1}^{\infty}  
  \frac{b_n^q}{n+m } \left( \frac{n}{m}\right)^{q \lambda} \right)^{1/q} \\ 
  & = \left( \sum_{m=1}^{\infty}  a_m^p \sum_{n=1}^{\infty} 
 \frac{m^{p \lambda} }{ (m+n) n^{p \lambda}}  \right)^{1/p}
 \left( \sum_{n=1}^{\infty}  b_n^q \sum_{m=1}^{\infty}  
\frac{n^{q \lambda} }{ (n+m) m^{q \lambda}}  \right)^{1/q}.
\end{align*}
We estimate the inner sums as follows.
For all $m,n \in \N$, the functions 
\[
f(x) =  \frac{m^{p \lambda} }{ (m+x) x^{p \lambda}} 
\qqand 
g(x) =  \frac{n^{q \lambda} }{ (n+x) x^{q \lambda}}. 
\]
are positive and decreasing for $x > 0$.  The integrals  
\[
 \int_0^{\infty} f(x)dx  \qqand \int_0^{\infty} g(x)dx  
\]
converge for $0 < \lambda <  1/p$ and for $0 < \lambda <  1/q$, respectively.  
With the substitution $x=mt$, we obtain 
\begin{align*}
 \sum_{n=1}^{\infty}  \frac{m^{p \lambda} }{ (m+n) n^{p \lambda}}   
 & =  \sum_{n=1}^{\infty} f(n) \leq \int_0^{\infty} f(x)dx \\
& =  \int_0^{\infty}  \frac{m^{p \lambda} }{ (m+x) x^{p \lambda}}  dx \\
&  =  \int_0^{\infty}  \frac{1}{ (1+t) t^{p \lambda}}  dt \\
& = \frac{\pi}{\sin\left(\pi p\lambda \right)}.
\end{align*}
With the substitution $x=nt$, we obtain 
\begin{align*}
\sum_{m=1}^{\infty}  \frac{n^{q \lambda} }{ (n+m) m^{q \lambda}} 
 & =  \sum_{m =1}^{\infty} g(m) \leq  \int_0^{\infty} g(x)dx \\ 
 & =  \int_0^{\infty}  \frac{n^{q \lambda} }{ (n+x) x^{q \lambda}} dx \\
& =  \int_0^{\infty}  \frac{1}{ (1+t) t^{q \lambda}}  dt \\
& = \frac{\pi}{\sin\left(\pi q\lambda \right)}. 
\end{align*}
Choosing
\[
\lambda = \frac{1}{pq}
\]
we obtain 
\[
\sin\left(\pi q\lambda \right) = \sin\left(\frac{\pi}{p}  \right) = \sin\left(\frac{\pi}{q}  \right) 
= \sin\left(\pi p\lambda \right)  
\]
and 
\begin{align*}
 \sum_{m=1}^{\infty} \sum_{n=1}^{\infty} \frac{a_m b_n}{m+n}  &
 \leq \left(     \sum_{m=1}^{\infty}  a_m^p  \frac{\pi}{\sin\left(\frac{\pi}{q}  \right)}   \right)^{1/p}
 \left( \sum_{n=1}^{\infty}  b_n^q   \frac{\pi}{\sin\left( \frac{\pi}{p}  \right)} \right)^{1/q} \\ 
& =  \frac{\pi} { \sin\left( \frac{\pi}{p}  \right) } 
\left( \sum_{m=1}^{\infty}  a_m^p \right)^{1/p}  \left( \sum_{n=1}^{\infty}  b_n^q \right)^{1/q}. 
\end{align*}
If $p = q = 2$, then 
\[
 \frac{\pi}{  \sin  \left(  \frac{\pi}{2} \right)  }  = \pi 
\]
gives Schur's best possible upper bound for Hilbert's inequality.  
This completes the proof. 
\end{proof}

Is there a converse to Hilbert's inequality that is similar to Landau's theorem? 
 For example, let $\mba = (a_m)_{m=1}^{\infty}$ be a sequence of real numbers such that, 
for all $\mbb = (b_n)_{n=1}^{\infty} \in \ell^q$,  the infinite series 
\[
\sum_{m=1}^{\infty} \sum_{n=1}^{\infty} \frac{a_m b_n}{m+n}
\]
converges.  Does this imply $\mba \in \ell^p$?



\begin{thebibliography}{30}




\bibitem{albu17}
N. Albuquerque, G. Ara\'{u}jo, D.  Pellegrino, and J. B. Seoane-Sep\'{u}lveda, 
H\"{o}lder's inequality: some recent and unexpected applications,
Bull. Belg. Math. Soc. Simon Stevin 24 (2017), 199--225. 




@\bibitem{borw08}
J. M. Borwein, Hilbert's inequality and {W}itten's zeta-function,
Amer. Math. Monthly 115 (2008), 125--137. 





\bibitem{chen-yang15}
Q. Chen and B. Yang,  
A survey on the study of {H}ilbert-type inequalities, 
Journal of Inequalities and Applications (2015). 



\bibitem{hard-litt-poly88}
G. H. Hardy,  J. E.  Littlewood, and G. P{\'o}lya, 
\emph{Inequalities}, Cambridge University Press, Cambridge, 1988. 

\bibitem{hell-topl06}
E. Hellinger and O. T\" oplitz, 
Grundlagen f\" ur eine Theorie der unendlichen Matrizen,
Nachrichten von der Gesellschaft der Wissenschaften zu G{\" o}ttingen, Mathematisch-Physikalische Klasse, (1906), 351--355.

\bibitem{knop56}
K. Knopp, \emph{Infinite Sequences and Series},
Dover Publications, New York, 1956.


\bibitem{land07}
E. Landau, {\" U}ber einen Konvergensatz, 
Nachrichten von der Gesellschaft der Wissenschaften zu G{\"o}ttingen, Mathematisch-Physikalische Klasse (1907), 25--27.



\bibitem{mali98}
L. Maligranda, 
Why {H}\"{o}lder's inequality should be called {R}ogers' inequality,
Math. Inequal. Appl. 1 (1998), 69--83.

\bibitem{mont-vaug73}
H. L. Montgomery and R. C. Vaughan, Hilbert's inequality, J. London Math. Soc.  (1973), 73--82.



\bibitem{ries13}
F. Riesz, \emph{Les syst{\` e}mes d'{\' e}quations lin{\'  e}aires {\` a} 
une infinit{\' e} d'inconnues}, Gauthier-Villars, Paris, 1913.



\bibitem{stee13}
J. M. Steele, 
\emph{The Cauchy-Schwarz Master Class}, Cambridge Univ. Press, 
Cambridge, 2013, pp. 155--165. 







\end{thebibliography}
\end{document}